\documentclass[12pt]{article}
\usepackage{amsmath,amssymb,amscd,verbatim,amsthm}
\usepackage{latexsym, epsfig, graphpap, graphics, mathrsfs}
\usepackage[varg]{pxfonts}

\begin{document}
 \noindent{\bf{Measures of noncompactness in some new lacunary difference sequence spaces}}

\vskip 0.5 cm

\noindent Ekrem Sava$\c{s}$

\noindent Istanbul Commerce University, 34840 Istanbul, Turkey

\noindent E-mail : ekremsavas@yahoo.com
\vskip 0.5 cm

\noindent Stuti Borgohain $^*${\footnote {The work of the
authors was carried under the Post Doctoral Fellow under National Board of Higher 
Mathematics, DAE, project No. NBHM/PDF.50/2011/64}}

\noindent Department of Mathematics

\noindent Indian Institute of Technology, Bombay

\noindent Powai:400076, Mumbai, Maharashtra; INDIA.
\vskip 0.5 cm

\noindent E-mail:stutiborgohain$@$yahoo.com
\vskip 1 cm

\noindent{\footnotesize {\bf{Abstract:}} In this research article, we establish some identities and estimates for the operator norms and the Hausdorff measures of noncompactness of certain operators on some lacunary difference sequence spaces defined by Orlicz function. Moreover, we apply our results to characterize some classes of compact operators on those spaces by using the Hausdorff measure of noncompactness.  \\


\vskip 0.3 cm

\noindent{\bf{AMS Classification No:} 40A05; 40A25; 40A30; 46B15; 46B45; 46B50.
40C05.}}

\vskip 1 cm

\section{Introduction}

Measures of noncompactness were first introduced and later on applied  in fixed point theory by Kuratowski [11] and Darbo [8]. Hausdorff measure of noncompactness was introduced by Goldenstein et. al and later on it was studied in broad sense by Eberhard Malkowsky et al. [5],  Feyzi Basar et al. [7], Eberhard Malkowsky and Ekrem Savas [6], Mohammed Mursaleen et al. [12,13] and many others. Some identities or estimations for the operator norms and  the Hausdorff measures of noncompactness of certain matrix operators on some sequence spaces were studied and estalished.\\

An important application of the Hausdorff measure of noncompactness of bounded linear operators between Banach spaces is the characterization of compact matrix transformations between $BK$ spaces. W.L.C. Sargent  proved that the characterizations of compact matrix operators between the classical sequence spaces in almost all cases. \\

Let $S$ and $M$ be subsets of a metric space $(X,d)$ and if for $\varepsilon > 0$ and for every $x \in M$ there exists $s \in S$ such that $d(x,s)<\varepsilon$, then $S$ is called an $\varepsilon$-net of $M$ in $X$. \\

Let $\mathscr{M}_X$  be a collection of all bounded subsets of a metric space $(X,d)$. The Hausdorff measure of non compactness of the set $Q$, denoted by $\chi(Q)$, is defined by, $\chi(Q)={\mbox{inf}}\{\varepsilon > 0: Q \mbox{~ has a finite ~} \varepsilon-\mbox{net in ~} X \}, \mbox{~where~} Q \in \mathscr{M}_X. $ The function $\chi: \mathscr{M}_X \rightarrow [0,\infty)$ is called the Hausdorff measure of noncompactness.\\

If $Q,Q_1$ and $Q_2$ are bounded subsets of a metric space $(X,d)$, then [see Malkowsky [5]),
$$\chi(Q)=0 \mbox{~if and only if~} Q \mbox{~is totally bounded~},$$
$$Q_1 \subset Q_2 \mbox{~implies ~} \chi(Q_1) \leq \chi(Q_2).$$ 

Further, the function $\chi$ has some additional properties connected with the linear structure, e.g.
$$\chi(Q_1+Q_2) \leq \chi(Q_1)+\chi(Q_2),$$
$$\chi(\alpha Q)= \vert \alpha \vert \chi(Q),\mbox{~for all ~} \alpha \in C.$$

Let $X$ and $Y$ be Banach spaces and $L \in B(X,Y)$. Then, the Hausdorff  measure of noncompactness of $L$, denoted by $\Vert L \Vert_\chi$, can be defined by,
\begin{equation}
\Vert L\Vert_\chi = \chi(L(S_X))= \chi(L(\overline B_X))
\end{equation}

and we have,

\begin{equation}
L \mbox{~is compact if and only if~} \Vert L\Vert_\chi=0.
\end{equation}

\section{Some preliminary concepts}

By a  lacunary sequence, we mean an increasing integer sequence $\theta = (k_r)$ such that $k_0=0$ and $h_r=k_r-k_{r-1} \rightarrow \infty$ as $r \rightarrow \infty$, where the intervals determined by $\theta$ will be denoted by $J_r=(k_{r-1}, k_r]$ and the ratio $\frac{k_r}{k_{r-1}}$ is defined by $\phi_r$. \\

For any lacunary sequence $\theta= (k_r)$, the space $N_\theta$ is defined as, (Freedman et al. [2])
$$N_\theta=\left\{(x_k): \displaystyle\lim_{r \rightarrow \infty} h_r^{-1} \displaystyle\sum_{k \in J_r} \vert x_k -L \vert =0, \mbox{~for some~} L \right\}.$$

The space $N_\theta$ is a $BK$ space with the norm,
$$\Vert (x_k) \Vert_\theta= \displaystyle\sup_r h_r^{-1} \displaystyle\sum_{k \in J_r} \vert x_k \vert.$$

An Orlicz function is defined as a function $M: [0,\infty )\rightarrow [0,\infty )$, which is continuous, non-decreasing and convex with
$M(0) = 0, M(x)>0$, for $x>0$ and $M(x)\rightarrow \infty$, as $x \rightarrow \infty$.\\

The idea of Orlicz function was used to construct the sequence space, [see Lindenstrauss and Tzafriri [10]
$$\ell_M=\left\{ (x_k) \in w: \displaystyle\sum_{k=1}^\infty M \left(\frac{\vert x_k \vert}{\rho}\right) < \infty, \mbox{~for some~} \rho>0 \right\}$$

which is a Banach space with the norm, called as Orlicz sequence space,
$$ \Vert x \Vert = \mbox{inf}\left\{ \rho>0: \displaystyle\sum_{k=1}^\infty M \left(\frac{\vert x_k \vert}{\rho}\right) \leq 1 \right\}.$$

The difference sequence spaces $\ell_\infty(\Delta), c(\Delta)$ and $c_0(\Delta)$ of crisp sets are defined as $Z(\Delta) = \{x = (x_k): (\Delta x_k) \in Z\},$ for $Z = \ell_\infty, c$ and $c_0$, where  $\Delta x = (\Delta x_k) = (x_k-x_{k+1})$, for all $k \in \mathbb{N}$, which can be a Banach space with $\Vert x \Vert_{\Delta}=\vert x_1 \vert + \displaystyle \sup_k \vert \Delta x_k \vert.$ \\

The generalized difference sequence spaces are defined as, for $m \geq 1$ and $n \geq 1$,

$$Z(\Delta_m^n) = \{x = (x_k): (\Delta_m^n x_k) \in Z\},\mbox{~for~} Z = \ell_\infty , c \mbox{~and~} c_0.$$

Let $X$ is any subset of $w$, then a matrix domain of an infinite matrix $A$ in $X$ is defined by, $X_A=\{x \in w: Ax \in X \}.$ If $x \supset \phi$ is a $BK$-space and $a=(a_k)\in w$, then we define, 

\begin{equation}
\Vert a \Vert_X^\ast=\displaystyle \sup_{x \in S_X} \left \vert \displaystyle \sum_{k=0}^\infty a_k x_k \right \vert.
\end{equation}

\section{The difference sequence spaces $c_0^\lambda(M, \Delta,s, \theta)$, $c^\lambda(M, \Delta,s, \theta)$ and  $\ell_\infty^\lambda(M, \Delta,s, \theta)$} 

Consider $\lambda=(\lambda_k)_{k=0}^\infty$ to be a strictly increasing sequence of positive reals such that $\lambda_k \rightarrow \infty$ as $k \rightarrow \infty$. We define the infinite matrix $\overline \Lambda =(\overline \lambda_{nk})_{n,k=0}^\infty$ by,

\begin{equation}
\overline \lambda_{nk}= \left\{
\begin{array}{l l}
\frac{(\lambda_k-\lambda_{k-1})-(\lambda_{k+1}-\lambda_k)}{\lambda_n}; ~~(k<n),\\
\frac{\lambda_n-\lambda_{n-1}}{\lambda_n};~~~~~~~~~~~~~~~(k=n),\\
0;~~~~~~~~~~~~~~~~~~~~~~~(k>n),
\end{array} \right. 
\end{equation}

where, we shall use the convention that any term with a negative subscript is equal to zero. Mursaleen and Noman [13] introduced the difference sequence spaces $c_0^\lambda(\Delta)$ and $\ell_\infty^\lambda(\Delta)$ as the matrix domains of the triangle $\overline\Lambda$ in the spaces $c_0$ and $\ell_\infty$ respectively.\\

In this paper, we study  the sequence spaces $c_0^\lambda(M, \Delta,s, \theta)$, $c^\lambda(M, \Delta,s, \theta)$ and  $\ell_\infty^\lambda(M, \Delta,s, \theta)$ and try to estimate for the operator norms and the Hausdorff measures of noncompactness of certain operators on these spaces. The spaces $c_0^\lambda(M, \Delta,s, \theta)$, $c^\lambda(M, \Delta,s, \theta)$ and  $\ell_\infty^\lambda(M, \Delta,s, \theta)$ are BK-spaces with the norm given by,

\begin{equation}
\Vert x \Vert_{\ell_\infty^\lambda(M,\Delta, s, \theta)}=\Vert \overline \Lambda (x) \Vert_{\ell_\infty(M, s, \theta)}=\mbox{inf}\left\{\rho>0:\displaystyle\lim_r \frac{1}{h_r} \displaystyle\sum_{k=1}^\infty \left(M \left(\frac{\vert \overline \Lambda_k(x) \vert}{\rho}\right)\right)^{s_k} \leq 1 \right\}.
\end{equation}

The $\beta-$duals of a subset $X$ of $w$ are respectively defined by,
$$X^\beta=\{a=(a_k) \in w: ax=(a_k x_k) \in cs,\mbox{~for all~} x=(x_k) \in X \}$$

{\bf{Lemma 3.1.}} Let $X$ denote any of the spaces $c_0^\lambda(M,\Delta,s , \theta)$ or $\ell_\infty^\lambda(M,\Delta, s, \theta)$. Then, we have,

\begin{equation}
\Vert a \Vert_X^\ast = \Vert \overline a \Vert_{\ell_1} = \displaystyle\sum_{k=0}^\infty \vert \overline a_k \vert< \infty
\end{equation}

for all $a = (a_k) \in X^\beta$, where,

\begin{equation}
\overline a_k=\lambda_k \left[\frac{a_k}{\lambda_k - \lambda_{k-1}} + \left(\frac{1}{\lambda_k- \lambda_{k-1}}-\frac{1}{\lambda_{k+1} -\lambda_k} \right) \displaystyle\sum_{j=k+1}^\infty a_j \right]; (k \in N)
\end{equation}

{\bf{Proof:}} Let $Y$ be the respective one of the spaces $c_0$ or $\ell_\infty$ .\\

Assume  $a=(a_k) \in X^\beta$ and  $y=\overline \Lambda (x)$ be the associated sequence defined by,

\begin{equation}
y_k=\displaystyle\sum_{j=0}^k \left(\frac{\lambda_j - \lambda_{j-1}}{\lambda_k} \right) (x_j -x_{j-1}); (k \in N).
\end{equation} 

Taking $y=\overline \Lambda (x)$ as the associated sequence, we have $\overline a = (\overline a_k) \in \ell_1$ such that for every $x=(x_k) \in X$, 

\begin{equation}
\displaystyle\sum_{k=0}^\infty a_k x_k = \displaystyle\sum_{k=0}^\infty \overline a_k y_k,
\end{equation}

Since $x \in S_X$ if and only if $y \in S_Y$, (followed by (7)), we can derive that, (by (1) and (9))\\

$\Vert a \Vert_X^\ast = \mbox{inf} \left \{\rho>0: \displaystyle\lim_r \frac{1}{h_r} \displaystyle\sum_{k=0}^\infty \left( M \left(\frac{\vert a_k x_k \vert }{\rho} \right)\right)^{s_k} \leq 1, x \in S_X \right\}$\\

$=\mbox{inf} \left\{ \rho>0: \displaystyle\lim_r \frac{1}{h_r}\displaystyle\sum_{k=0}^\infty \left(M \left( \frac{\vert \overline a_k y_k \vert }{\rho} \right)\right)^{s_k} \leq 1, y \in S_Y \right\}$\\

$=\Vert \overline a \Vert_Y^\ast$\\

It is known that  $\Vert . \Vert _X^\ast = \Vert . \Vert _{X^\beta}$ on $X^\beta$, where $\Vert . \Vert_{X^\beta}$  denotes the natural norm on the dual space $X^\beta$ and $X= c_0, c,\ell_\infty$ or $\ell_p(1 \leq p<\infty)$.\\
 
So if $\overline a \in \ell_1$, we obtain that $\Vert a \Vert_X^\ast = \Vert \overline a \Vert_Y^\ast = \Vert\overline a \Vert_{\ell_1}< \infty$, which concludes the proof.\\

Let $A=(a_{nk})$ be an infinite matrix and $\overline A =(\overline a_{nk})$ is the associated matrix defined by,

\begin{equation}
\overline a_{nk} = \lambda_k \left[ \frac{a_{nk}}{\lambda_k - \lambda_{k-1}}+\left(\frac{1}{\lambda_k - \lambda_{k-1}}- \frac{1}{\lambda_{k+1} -\lambda_k}\right) \displaystyle\sum_{j=k+1}^\infty a_{nj} \right]; (n,k \in N)
\end{equation}
\\
{\bf{Lemma 3.2.}} Let $X$ be any of the spaces $c_0^\lambda(M,\Delta,s,\theta)$ or $\ell_\infty^\lambda(M,\Delta)$and $Z$ be a  sequence space.. If $A \in (X,Z)$, then $\overline A \in (Y,Z)$ such that $Ax =\overline A y$ for all sequences $x \in X$ and $y \in Y$, here $Y$ is the respective one of the spaces $c_0$ or $\ell_\infty$.\\

{\bf{Proof.}} Suppose that $A \in (X,Z)$.\\

For any sequence  $x=(x_k) \in w$ and the associated sequence $y=\overline \Lambda(x)$ defined in (8), we have,

\begin{equation}
x_k=\displaystyle\sum_{j=0}^k \left(\frac{\lambda_j y_j - \lambda_{j-1} y_{j-1}}{\lambda_j - \lambda_{j-1}} \right) ; (k \in N).
\end{equation}

Then, $A_n \in X^\beta$ for all $n \in N$. Also, $\overline A_n \in \ell_1=Y^\beta$ for all $n \in N$ and the equality $Ax=\overline A y$, followed by equations (8), (9) and (10).  Hence, $\overline A y \in Z$. Further, by (11), we get that every $y \in Y$ is the associated sequence of some $x \in X$. Thus, it can be deduced that $\overline A \in (Y,Z)$, which completes the proof.\\

{\bf{Lemma 3.3.}} Let $X$ be any of the spaces $c_0^\lambda(M, \Delta, s, \theta)$ or $\ell_\infty^\lambda(M,\Delta, s, \theta)$, $A=(a_{nk})$ an infinite matrix and $\overline A=(\overline a_{nk})$ the associated marix. If $A$ is in any of the classes $(X,c_0),(X,c)$ or $(X,\ell_\infty)$, then,

$$\Vert L_A \Vert=\Vert A \Vert_{(X,\ell_\infty)}=\displaystyle\sup_n\left\{\mbox{inf}\left(\rho>0:\displaystyle\lim_r \frac{1}{h_r}\displaystyle\sum_{k=0}^\infty\left( M\left(\frac{\vert \overline a_{nk} \vert }{\rho}\right)\right)^{s_k} \leq 1 \right)\right\}< \infty$$

\section{Compact operators on the spaces $c_0^\lambda(M, \Delta, s, \theta)$ and $\ell_\infty^\lambda(M, \Delta,s, \theta)$}

In this section, we are trying to establish some identities or estimates for the Hausdorff measures of noncompactness of certain matrix operators on the spaces $c_0^\lambda(M, \Delta,s, \theta)$ and $\ell_\infty^\lambda(M, \Delta,s,\theta)$. Also, the results obtained by examining these sequence spaces are applied to characterize some classes of compact operators on those spaces. \\

{\bf{Remark: 4.1.}} Let $X$ denote any of the spaces $c_0$ or $\ell_\infty$. If $A \in (X,c)$, then we have,

\begin{itemize}

\item $\alpha_k=\displaystyle\lim_{n \rightarrow\infty} a_{nk}~\mbox{exists for every~} k \in N,$

\item $\alpha=(\alpha_k) \in \ell_1,$

\item $\displaystyle\sup_n \left(\displaystyle\sum_{k=0}^\infty \vert a_{nk}-\alpha_k \vert \right) < \infty,$

\item $\displaystyle\lim_{n \rightarrow \infty} A_n(x)=\displaystyle\sum_{k=0}^\infty \alpha_k x_k \mbox{~for all~} x=(x_k)\in X.$

\end{itemize}   

{\bf{Theorem 4.2.}} Assume $A=(a_{nk})$ be an infinite matrix and $\overline A = (\overline a_{nk})$ the associated matrix defined by (10). Then, we have the following results on the Hausdorff measures of noncompactness on the sequence spaces $X=c_0^\lambda(M,\Delta,s,\theta),c^\lambda(M,\Delta,s,\theta)$ or $\ell_\infty^\lambda(M,\Delta,s,\theta)$.\\

(a) If $A \in (X,c_0)$, then,

\begin{equation}
\Vert L_A \Vert_\chi = \displaystyle\lim_{n \rightarrow \infty} \sup \left(\mbox{inf}\left\{\rho>0:\displaystyle\lim_r \frac{1}{h_r}\displaystyle\sum_{k =0}^\infty \left(M\left(\frac{\vert\overline a_{nk} \vert}{\rho} \right)\right)^{s_k} \leq 1 \right\}\right).
\end{equation}

(b) If $A \in (X,c)$, then,\\

$\frac{1}{2} \displaystyle\lim_{n \rightarrow \infty}\sup \left(\mbox{inf}\left\{\rho>0: \displaystyle\lim_r \frac{1}{h_r}\displaystyle\sum_{k=0}^\infty \left(M \left(\frac{\vert \overline a_{nk}-\overline \alpha_k \vert}{\rho}\right)\right)^{s_k} \leq 1 \right\}\right)$

\begin{equation}
\leq \Vert L_A \Vert_\chi \leq  \displaystyle\lim_{n \rightarrow \infty}\sup \left( \mbox{inf}\left\{\rho>0: \displaystyle\lim_r \frac{1}{h_r}\displaystyle\sum_{k=0}^\infty \left(M \left(\frac{\vert \overline a_{nk}-\overline \alpha_k \vert}{\rho}\right)\right)^{s_k} \leq 1 \right\}\right),
\end{equation}   

where $\overline \alpha_k = \displaystyle\lim_{n \rightarrow \infty} \overline a_{nk}$ for all $k \in N$.\\

(c) If $A \in (X,\ell_\infty)$, then,

\begin{equation}
0 \leq \Vert L_A \Vert_\chi \leq \displaystyle\lim_{n \rightarrow \infty} \sup \left(\mbox{inf} \left\{\rho>0: \displaystyle\lim_r \frac{1}{h_r}\displaystyle\sum_{k=0}^\infty \left(M\left(\frac{\vert \overline a_{nk}\vert}{\rho} \right)\right)^{s_k} \leq 1 \right\} \right).
\end{equation}

{\bf{Proof:}} Following Lemma 3.3, it can be easily proved that the expressions in (12) and (13) exist. \\

Similarly, following Remark 4.1. and Lemma 3.2, we can deduce that the expression in (15) also  exists.\\

Let  $X \supset \phi$ and $Y$ be $BK$-spaces. Then, we have, $(X,Y) \subset B(X,Y)$, that is, every matrix $A \in (X,Y)$ defines an operator $L_A \in B(X,Y)$ by $L_A(x)=Ax$ for all $x \in X$.\\

So,
\begin{equation}
 \Vert L_A \Vert_\chi = \chi(AS), \mbox{~where  $S=S_X$}, \mbox{~ by(1)}
\end{equation}
 
Let $P_r:c_0 \rightarrow c_0 ~(r \in N)$ be the operator defined by $P_r(x)=(x_0,x_1,...x_r,0,0,..)$ for all $x=(x_k) \in c_0$. Then, we have, for $Q \in \mathscr{M}_{c_0}$, 

$$\chi(Q)=\displaystyle\lim_{r \rightarrow \infty}(\displaystyle\sup_{x \in Q} \Vert(I-P_r)(x) \Vert_{\ell_\infty}).$$

where $I$ is the identity operator on $c_0$.

Further, every $z=(z_n) \in c$ has a unique representation as $z=\overline z e + \displaystyle\sum_{n=0}^\infty (z_n-\overline z) e^{(n)}$, where $\overline z = \displaystyle\lim_{n \rightarrow \infty} z_n$. The projectors $P_r:c \rightarrow c~(r \in N)$ are obtained by,
\begin{equation}
 P_r(z)=\overline z e +\displaystyle\sum_{n=0}^r (z_n-\overline z) e^{(n)}; ~(r \in N)
\end{equation}

for all $z=(z_n) \in c$ with $\overline z= \displaystyle\lim_{n \rightarrow \infty} z_n$. \\

Let $AS \in \mathscr{M}_{c_0}$. Then, 
\begin{equation}
 \chi(AS)=\displaystyle\lim_{r \rightarrow \infty} (\sup_{x \in S} \Vert (I-P_r)(A_x) \Vert_{\ell_\infty}),
\end{equation}

where $P_r:c_0 \rightarrow c_0 (r \in N)$.\\

This implies that,

\begin{equation}
\Vert (I-P_r)(Ax) \Vert_{\ell_\infty} = \displaystyle\sup_{n > r} \vert A_n(x) \vert, \mbox{~ for all~}x \in X \mbox{~and every~} r \in N.
\end{equation}

For an infinite matrix $A=(a_{nk})_{n,k=0}^\infty$, we have the $A$-transform of $x$ as the sequence $Ax=(A_n(x))_{n=0}^\infty$, where $A_n(x)=\displaystyle\sum_{k=0}^\infty a_{nk} x_k$, for $x \in w$ and $n \in N$ .\\

Thus, we get, (by (3) and Lemma 3.1)

$$\displaystyle\sup_{x \in S} \Vert (I-P_r)(Ax) \Vert_{\ell_\infty}=\displaystyle\sup_{n>r} \Vert A_n \Vert_{(\ell_\infty^\lambda(M,\Delta,s,\theta))}^\ast = \displaystyle\sup_{n>r} \Vert \overline A_n \Vert_{\ell_1}, {\mbox{~ for every~}} r \in N.$$

Which implies that, (using above with(17))

$$\chi(AS)=\displaystyle\lim_{r \rightarrow \infty}(\displaystyle\sup_{n>r} \Vert \overline A_n \Vert_{\ell_1} ) = \displaystyle\lim_{n \rightarrow \infty}\sup \Vert \overline A_n \Vert_{\ell_1}.$$

This concludes the proof of (a). \\

To prove (b), let us take $Q \in \mathscr{M}_c$ and $P_r:c \rightarrow c ~(r \in N)$ be the projector onto the linear span of $\{e,e^{(0)},e^{(1)},...e^{(r)}\}$. Then, we have

$$\frac{1}{2}. \displaystyle\lim_{r \rightarrow \infty} \left(\displaystyle\sup_{x \in Q} \Vert (I-P_r)(x) \Vert_{\ell_\infty} \right) \leq \chi(Q) \leq \displaystyle\lim_{r \rightarrow \infty} \left(\displaystyle\sup_{x \in Q} \Vert (I-P_r)(x) \Vert_{\ell_\infty} \right),$$

where $I$ is the identity operator on $c$.\\

Since we have $AS \in \mathscr{M}_c$. We can get an estimate for the value of $\chi(AS)$ in (15). For this, let $P_r: c \rightarrow c(r \in N)$ be the projectors defined by (16).\\

Then, we have for every $r \in N$ that,

$$(I-P_r)(z)=\displaystyle\sum_{n=r+1}^\infty(z_n-\overline z)e^{(n)} $$

and hence,

\begin{equation}
 \Vert (I-P_r)(z)\Vert_{\ell_\infty} = \displaystyle\sup_{n > r} \vert z_n -\overline z \vert
\end{equation}
 
for all $z=(z_n) \in c$ and every $r \in N$, where $\overline z = \displaystyle\lim_{n \rightarrow \infty} z_n$ and $I$ is the identity operator on $c$.\\

 By using (15), we obtain,

\begin{equation}
 \frac{1}{2} \displaystyle\lim_{r \rightarrow \infty} (\displaystyle\sup_{x \in S} \Vert (I-P_r)(Ax) \Vert_{\ell_\infty}) \leq \Vert L_A \Vert_\chi \leq \displaystyle\lim_{r \rightarrow \infty} (\displaystyle\sup_{x \in S} \Vert (I-P_r)(Ax) \Vert_{\ell_\infty})
\end{equation}

On the other hand, it is given that $X=c_0^\lambda(M, \Delta,s,\theta)$ or $X=\ell_\infty^\lambda(M,\Delta,s,\theta)$, and let $Y$ be the respective one of the spaces $c_0$ or $\ell_\infty$. Also,  let $y \in Y$ be the associated sequence defined by (8). Since $A \in (X,c)$, we have from Lemma 3.2 that $\overline A \in (Y,c)$ and $Ax=\overline A y$. Further,  we have the limits $\overline \alpha_k=\displaystyle\lim_{n \rightarrow \infty} \overline a_{nk}$ exist for all $k, \overline \alpha = (\overline \alpha_k) \in \ell_1=Y^\beta$ and $\displaystyle\lim_{n \rightarrow \infty} \overline A_n(y)=\displaystyle\sum_{k=0}^\infty \overline \alpha_k y_k$. (Remark 4.1) \\

Consequently,\\

$\Vert (I-P_r)(Ax) \Vert_{\ell_\infty} = \Vert (I-P_r)(\overline A y) \Vert_{\ell_\infty}$\\

$~~~~=\displaystyle\sup_{n>r} \vert \overline A_n(y)-\displaystyle\sum_{k=0}^\infty \overline \alpha_k y_k \vert$\\

$~~~~=\displaystyle\sup_{n>r} \vert \displaystyle\sum_{k=0}^\infty (\overline a_{nk}-\overline \alpha_k)y_k \vert$, for all $r \in N$ (by (21)).\\

Moreover, since $x \in S=S_X$ if and only if $y \in S_Y$, we get,\\

$\displaystyle\sup_{x \in S} \Vert (I-P_r)(Ax) \Vert_{\ell_\infty}=\displaystyle\sup_{n>r}\left(\displaystyle\sup_{y \in S_Y} \vert \displaystyle\sum_{k=0}^\infty (\overline a_{nk}-\overline \alpha_k)y_k \vert\right)$\\

$~~~~=\displaystyle\sup_{n>r} \Vert \overline A_n \overline \alpha \Vert_Y^\ast$\\

$~~~~=\displaystyle\sup_{n>r} \Vert \overline A_n \overline \alpha \Vert_{\ell_1}$\\

for all $r \in N$. \\

This concludes the proof.

Finally, to prove (c), let us define the operators $P_r:\ell_\infty \rightarrow \ell_\infty (r \in N)$ as in the proof of part (a) for all $x=(x_k) \in \ell_\infty$. Then, we have,

$$AS \subset P_r(AS) + (I-P_r)(AS); (r \in N).$$

Thus, following the elementary properties of the function $\chi$, we have, \\

$0 \leq \chi(AS) \leq \chi(P_r(AS)) + \chi((I-P_r)(AS))$\\

$~~~=\chi((I-P_r)(AS))$\\

$~~~\leq \displaystyle\sup_{x \in S} \Vert (I-P_r)(Ax) \Vert_{\ell_\infty}$\\

$=\displaystyle\sup_{n>r} \Vert \overline A_n \Vert_{\ell_1}$, for all $r \in N$.\\

Hence,\\

$0 \leq \chi(AS) \leq \displaystyle\lim_{r \rightarrow \infty}(\displaystyle\sup_{n>r} \Vert \overline A_n \Vert_{\ell_1})$\\

$=\displaystyle\lim_{n \rightarrow \infty} \displaystyle\sup \Vert \overline A_n \Vert_{\ell_1}$.\\

Combining this together with (15), imply (14) ,which completes the proof.\\

{\bf{Corollary 4.3.}} Let $X$ denote any of the spaces $c_0^\lambda(M,\Delta,s,\theta)$ or $\ell_\infty^\lambda(M,\Delta,s,\theta)$. Then, we have,\\

(a) If $A \in (X,c_0)$, then,\\

$L_A$ is compact iff $\displaystyle\lim_{n \rightarrow \infty}\left(\mbox{inf} \left\{\rho>0: \displaystyle\lim_r \frac{1}{h_r}\displaystyle\sum_{k=0}^\infty M\left(\frac{\vert \overline a_{nk} \vert}{\rho}\right)^{s_k}\leq 1\right\}\right) =0$.\\

(b) If $A \in (X,c)$, then,\\

$L_A$ is compact iff $\displaystyle\lim_{n \rightarrow \infty} \left(\mbox{inf}\left\{\rho>0: \displaystyle\lim_r \frac{1}{h_r}\displaystyle\sum_{k=0}^\infty \left(M \left(\frac{\vert \overline a_{nk}-\overline \alpha_k \vert}{\rho}\right)\right)^{s_k}\leq 1 \right\}\right) =0$ where 

$\overline \alpha_k = \displaystyle\lim_{n \rightarrow \infty} \overline a_{nk}$ for all $k \in N$.\\

(c) If $A \in (X,\ell_\infty)$, then,\\

$L_A$ is compact if $\displaystyle\lim_{n \rightarrow \infty}\left(\mbox{inf}\left\{\rho>0: \displaystyle\lim_r \frac{1}{h_r}\displaystyle\sum_{k=0}^\infty \left(M\left(\frac{\vert \overline a_{nk} \vert}{\rho}\right)\right)^{s_k}\leq 1 \right\}\right) =0$.
  
\section{Some applications}

By applying the previous results, in this section, we are trying to establish some identities or estimates for the operator norms and the Hausdorff measure of non-compactness of certain matrix operators that map any of the spaces $c_0^\lambda (M,\Delta,s,\theta), c^\lambda(M,\Delta,s,\theta)$ and $\ell_\infty^\lambda(M,\Delta,s,\theta)$ into the matrix domains of  triangles in the spaces $c_0,c $ and $\ell_\infty$. Further, we deduce the necessary and sufficient conditions for such operators to be compact.\\

\noindent {\bf{Lemma 5.1.}} Let $T$ be a triangle. Then, we have,

\begin{itemize}

\item For arbitrary subsets $X$ and $Y$ of $w$, $A \in (X,Y_T)$ if and only if $B=TA \in (X,Y)$.

\item Further, if $X$ and $Y$ are $BK$ spaces and $A \in (X,Y_T)$, then $\Vert L_A \Vert=\Vert L_B \Vert$.

\end{itemize} 

Throughout, we assume that $A=(a_{nk})$ is an infinite matrix and $T=(t_{nk})$ is a triangle, and we define the matrix $B=(b_{nk})$ by $b_{nk}=\displaystyle\sum_{m=0}^n t_{nm} a_{mk}; (n,k \in N)$, that is $B=TA$ and hence,

$$B_n=\displaystyle\sum_{m=0}^n t_{nm}A_m = \left(\displaystyle\sum_{m=0}^n t_{nm} a_{mk}\right)_{k=0}^\infty; (n \in N).$$

Consider $\overline A=(\overline a_{nk})$ and $\overline B=(\overline b_{nk})$ be the associated matrices of $A$ and $B$, respectively. Then it can easily be seen that,

$$\overline b_{nk}=\displaystyle\sum_{m=0}^n t_{nm} \overline a_{mk};(n,k\in N).$$

Hence, $\overline B_n=\displaystyle\sum_{m=0}^n t_{nm} \overline A_m =\left(\displaystyle\sum_{m=0}^n t_{nm} \overline a_{mk}\right)_{k=0}^\infty; (n \in N)$.  \\

Moreover, we define the sequence $\overline a=(\overline a_k)_{k=0}^\infty$ by,

$$\overline a_k=\displaystyle\lim_{n \rightarrow \infty}\left(\displaystyle\sum_{m=0}^n t_{nm} \overline a_{mk} \right); (k \in N)$$

provided the above limits exist for all $k \in N$ which is the case whenever $A \in (c_0^\lambda(M,\Delta,s,\theta), c_T)$ or $A\in (\ell_\infty^\lambda(M,\Delta,s,\theta),c_T)$ by lemmas 5.1, 3.2 and Remark 4.1.\\ 

Now using the above results, we have the following results:\\

{\bf{Theorem 5.2.}} Let $X$ be any of the spaces $c_0^\lambda(M,\Delta,s,\theta)$ or $\ell_\infty^\lambda(M,\Delta,s,\theta)$, $T$ a traingle and $A$ an infinite matrix. If $A$ is in any of the classes $(X,(c_0)_T),(X,c_T)$ or $(X,(\ell_\infty)_T)$, then\\

$\Vert L_A \Vert =\Vert A \Vert_{(X,(\ell_\infty)_T)}=\displaystyle\sup_n \left\{ \mbox{inf}\left( \rho>0: \displaystyle\lim_r \frac{1}{h_r} \displaystyle\sum_{k=0}^\infty \left\vert \displaystyle \sum_{m=0}^n \left (M \left( \frac{\vert t_{nm} \overline a_{mk} \vert}{\rho}\right)\right)^{s_k} \right \vert \leq 1 \right)\right\} <\infty.$\\

{\bf{Theorem 5.3.}}  Let $T$ be a triangle. If either $A \in (\ell_\infty^\lambda(M,\Delta,s \theta),(c_0)_T)$ or $A \in (\ell_\infty^\lambda(M,\Delta,s,\theta), c_T)$ then $L_A$ is compact.\\

{\bf{Theorem 5.4.}}  Let $T$ be a triangle. Then, we have,

\begin{enumerate}

\item If $A \in (c_0^\lambda(M,\Delta,s,\theta),(c_0)_T)$, then,\\

$\Vert L_A \Vert_\chi = \lim \displaystyle\sup_{n \rightarrow \infty} \left\{ \mbox{inf} \left(\rho>0: \displaystyle\lim_r \frac{1}{h_r} \displaystyle\sum_{k=0}^\infty \left \vert \displaystyle\sum_{m=0}^n \left( M \left( \frac{\vert t_{nm} \overline a_{mk} \vert}{\rho} \right)\right) ^{s_k} \right\vert \leq 1 \right)\right\}.$\\

and $L_A$ is compact if and only if 

$$\displaystyle\lim_{n \rightarrow \infty} \left\{ \mbox{inf} \left(\rho>0: \displaystyle\lim_r \frac{1}{h_r} \displaystyle\sum_{k=0}^\infty \left\vert \displaystyle\sum_{m=0}^n \left( M \left(\frac{\vert t_{nm} \overline a_{mk} \vert}{\rho} \right)\right)^{s_k} \right \vert \leq 1 \right)\right\} =0.$$

\item If $A \in (c_0^\lambda(M,\Delta,s,\theta), c_T)$, then

$\frac{1}{2}. \displaystyle\lim\sup_{n \rightarrow \infty} \left\{ \mbox{inf} \left(\rho>0: \displaystyle\lim_r \frac{1}{h_r} \displaystyle\sum_{k=0}^\infty \left \vert \displaystyle\sum_{m=0}^n \left( M \left(\frac{\vert t_{nm} \overline a_{mk}-\overline a_k \vert}{\rho} \right)\right)^{s_k} \right \vert \leq 1 \right)\right\}$\\

~~$ \leq \Vert L_A \Vert_\chi \leq \displaystyle\lim\sup_{n \rightarrow \infty} \left\{ \mbox{inf} \left(\rho>0: \displaystyle\lim_r \frac{1}{h_r} \displaystyle\sum_{k=0}^\infty \left \vert \displaystyle\sum_{m=0}^n \left( M \left(\frac{\vert t_{nm} \overline a_{mk}-\overline a_k \vert}{\rho} \right)\right)^{s_k} \right \vert \leq 1 \right)\right\}$

and $L_A$ is compact if and only if
$$\displaystyle\lim_{n \rightarrow \infty} \left\{ \mbox{inf} \left(\rho>0: \displaystyle\lim_r \frac{1}{h_r} \displaystyle\sum_{k=0}^\infty \left \vert \displaystyle\sum_{m=0}^n \left( M \left(\frac{\vert t_{nm} \overline a_{mk} - \overline a_k \vert}{\rho} \right)\right)^{s_k} \right \vert \leq 1 \right)\right\}=0.$$ 

\item If either $A \in (c_0^\lambda(M,\Delta,s,\theta),(\ell_\infty)_T)$ or $A \in (\ell_\infty^\lambda(M,\Delta,s,\theta),(\ell_\infty)_T)$, then
$$0 \leq \Vert L_A \Vert_\chi \leq \displaystyle\lim\sup_{n \rightarrow \infty}\left\{ \mbox{inf} \left(\rho>0: \displaystyle\lim_r \frac{1}{h_r} \displaystyle\sum_{k=0}^\infty \left \vert \displaystyle\sum_{m=0}^n \left( M \left(\frac{\vert t_{nm} \overline a_{mk} \vert}{\rho} \right)\right)^{s_k} \right \vert \leq 1 \right)\right\} $$

and $L_A$ is compact if

$$\displaystyle\lim_{n \rightarrow \infty} \left\{ \mbox{inf} \left(\rho>0: \displaystyle\lim_r \frac{1}{h_r} \displaystyle\sum_{k=0}^\infty \left \vert \displaystyle\sum_{m=0}^n \left( M \left(\frac{\vert t_{nm} \overline a_{mk} \vert}{\rho} \right)\right)^{s_k} \right \vert \leq 1 \right)\right\} =0.$$
 
\end{enumerate}

\vskip 0.5 cm

Particular cases: Let $\lambda'=(\lambda'_k)_{k=0}^\infty$ be a strictly increasing sequence of positive reals tending to infinity and $\Lambda'=(\lambda'_{nk})$ be the triangle defined by (4), with the sequence $\lambda'$ instead of $\lambda$. Also, let $c_0^{\lambda'}(M,\Delta,s,\theta), c^{\lambda'}(M,\Delta,s,\theta)$ and $\ell_\infty^{\lambda'}(M,\Delta,s,\theta)$ be the matrix  domains of the triangle $\Lambda'$ in the spaces $c_0,c$ and $\ell_\infty$ respectively.\\

{\bf{Particular Case 5.5.}} Let $X$ be any of the spaces $c_0^\lambda(M,\Delta,s,\theta)$ or $\ell_\infty^\lambda(M,\Delta,s,\theta)$ and $A$ an infinite matrix. If $A$ is in any of the classes $(X,c_0^{\lambda'}(M,\Delta,s,\theta)), (X,c^{\lambda'}(M,\Delta,s,\theta))$ or $(X, \ell_\infty^{\lambda'}(M,\Delta,s,\theta))$, then

$$\Vert L_A \Vert = \Vert A \Vert_{(X,\ell_\infty^{\lambda'}(M,\Delta,s,\theta))}=\displaystyle\sup_n \left\{ \mbox{inf} \left(\rho>0: \displaystyle\lim_r \frac{1}{h_r} \displaystyle\sum_{k=0}^\infty \left\vert\displaystyle\sum_{m=0}^n \left( M \left(\frac{\vert \lambda'_{nm} \overline a_{mk} \vert}{\rho} \right)\right)^{s_k} \right\vert \leq 1 \right)\right\} .$$

{\bf{Particular Case 5.6.}} If either $A \in (\ell_\infty^\lambda(M,\Delta,s,\theta), c_0^{\lambda'}(M,\Delta,s,\theta))$ or $A \in (\ell_\infty^\lambda(M,\Delta,s,\theta), c^{\lambda'}(M,\Delta,s,\theta))$, then $L_A$ is compact. \\

Similarly, we get some identities or estimates for the Hausdorff measures of noncompactness of operators given by matrices in the classes $(c_0^\lambda(M,\Delta,s,\theta), c_0^{\lambda'}(M,\Delta,s,\theta)), (c_0^\lambda(M,\Delta,s,\theta), c^{\lambda'}(M,\Delta,s,\theta)), (c_0^\lambda(M,\Delta,s,\theta),\ell_\infty^{\lambda'}(M,\Delta,s,\theta))$ and $(\ell_\infty^\lambda(M,\Delta,s,\theta),\ell_\infty^{\lambda'}(M,\Delta,s,\theta))$,and deduce the necessary and suffucient (or only sufficient ) conditions for such operators to be compact.\\

Let $bs, cs$ and $cs_0$ be the spaces of all sequences associated with bounded, convergent and null series,respectively. Then, we have the following results associated with the these sequence spaces,\\

{\bf{Corollary 5.7.}} Let $X$ be any of the spaces $c_0^\lambda(M,\Delta,s,\theta)$ or $\ell_\infty^\lambda(M,\Delta,s,\theta)$ and $A$ an infinite matrix. If $A$ is in any of the classes $(X, cs_0), (X, cs)$ or $(X, bs)$, then,

$$\Vert L_A \Vert = \Vert A \Vert_{(X, bs)}=\displaystyle\sup_n \left\{ \mbox{inf} \left(\rho>0: \displaystyle\lim_r \frac{1}{h_r} \displaystyle\sum_{k=0}^\infty \left\vert \displaystyle\sum_{m=0}^n \left( M \left(\frac{\vert \overline a_{mk} \vert}{\rho} \right)\right)^{s_k} \right\vert \leq 1 \right)\right\} < \infty$$. \\

{\bf{Corollary 5.8.}} If either $A \in (\ell_\infty^\lambda(M,\Delta,s,\theta), cs_0)$ or $A \in (\ell_\infty^\lambda(M,\Delta,s,\theta), cs)$,then $L_A$ is compact.\\

{\bf{Corollary 5.9.}} We have

\begin{enumerate}

\item If $A \in (c_0^\lambda(M,\Delta,s,\theta), cs_0)$, then

$$\Vert L_A \Vert_\chi = \displaystyle\lim\sup_{n \rightarrow \infty} \left\{ \mbox{inf} \left(\rho>0: \displaystyle\lim_r \frac{1}{h_r} \displaystyle\sum_{k=0}^\infty \left \vert\displaystyle\sum_{m=0}^n \left( M \left(\frac{\vert \overline a_{mk} \vert}{\rho} \right)\right)^{s_k} \right \vert\leq 1 \right)\right\}$$

and $L_A$ is compact if and only if 

$$\displaystyle\lim_{n \rightarrow \infty} \left\{ \mbox{inf} \left(\rho>0: \displaystyle\lim_r \frac{1}{h_r} \displaystyle\sum_{k=0}^\infty \left\vert \displaystyle\sum_{m=0}^n \left( M \left(\frac{\vert \overline a_{mk} \vert}{\rho} \right)\right)^{s_k} \right\vert\leq 1 \right)\right\} =0.$$

\item If $A \in (c_0^\lambda(M,\Delta,s,\theta), cs)$, then

$\frac{1}{2}.\displaystyle\lim\sup_{n \rightarrow \infty} \left\{ \mbox{inf} \left(\rho>0: \displaystyle\lim_r \frac{1}{h_r} \displaystyle\sum_{k=0}^\infty \left\vert \displaystyle\sum_{m=0}^n \left( M \left(\frac{\vert \overline a_{mk} - \overline a_k \vert}{\rho} \right)\right)^{s_k} \right\vert \leq 1 \right)\right\}$\\

~~~~~$\leq \Vert L_A \Vert_\chi \leq  \displaystyle\lim\sup_{n \rightarrow \infty} \left\{ \mbox{inf} \left(\rho>0: \displaystyle\lim_r \frac{1}{h_r} \displaystyle\sum_{k=0}^\infty \left\vert \displaystyle\sum_{m=0}^n \left( M \left(\frac{\vert \overline a_{mk} - \overline a_k \vert}{\rho} \right)\right)^{s_k} \right\vert \leq 1 \right)\right\}$\\

and $L_A$ is compact if and only if 

$$\displaystyle\lim_{n \rightarrow \infty} \left\{ \mbox{inf} \left(\rho>0: \displaystyle\lim_r \frac{1}{h_r} \displaystyle\sum_{k=0}^\infty \left\vert \displaystyle\sum_{m=0}^n \left( M \left(\frac{\vert \overline a_{mk} - \overline a_k \vert}{\rho} \right)\right)^{s_k} \right\vert \leq 1 \right)\right\} =0,$$

where $\overline a_k = \displaystyle\lim_{n \rightarrow \infty} (\displaystyle\sum_{m=0}^n \overline a_{mk})$ for all $k\in N.$

\item If either $A \in (c_0^\lambda(M,\Delta,s,\theta), bs)$ or $A \in (\ell_\infty^\lambda(M,\Delta,s,\theta),bs)$, then

$$0 \leq \Vert L_A \Vert_\chi \leq \displaystyle\lim\sup_{n \rightarrow \infty} \left\{ \mbox{inf} \left(\rho>0: \displaystyle\lim_r \frac{1}{h_r} \displaystyle\sum_{k=0}^\infty \left\vert \displaystyle\sum_{m=0}^n \left( M \left(\frac{\vert \overline a_{mk} \vert}{\rho} \right)\right)^{s_k} \right\vert \leq 1 \right)\right\} $$

and $L_A$ is compact if 

$$\displaystyle\lim_{n \rightarrow \infty} \left\{ \mbox{inf} \left(\rho>0: \displaystyle\lim_r \frac{1}{h_r} \displaystyle\sum_{k=0}^\infty \left\vert \displaystyle\sum_{m=0}^n \left( M \left(\frac{\vert \overline a_{mk} \vert}{\rho} \right)\right)^{s_k} \right\vert \leq 1 \right)\right\}=0.$$

\end{enumerate}


\vskip 0.5 cm

\begin{thebibliography}{50}

\bibitem{article} A. Alotaibi, B.  Hazarika and  S. A. Mohiuddine (2014).  On lacunary statistical convergence of double sequences in locally solid Riesz spaces. {\it{J. Comput. Anal. Appl.}} 17 (1) : 156–165. 

\bibitem{article} A.R. Freedman, J.J.  Sember and  M. Raphael (1978).  Some Ces$\grave{a}$aro type summability space. {\it{Proceedings of the London Mathematical Society}}. 37(3):508–520.

\bibitem{article} A.  G\"{o}khan, M.  Et  and  M. Mursaleen (2009). Almost lacunary statistical and strongly almost lacunary convergence of sequences of fuzzy numbers. {\it{Mathematical and Computer Modelling}}.49:548-555.

\bibitem{article} B. Hazarika, S.A.  Mohiuddine  and M. Mursaleen (2014).  Some inclusion results for lacunary statistical convergence in locally solid Riesz spaces. {\it{Iran. J. Sci. Technol. Trans. A Sci.}}. 38(1):61–68. 

\bibitem{article} E. Malkowsky (2010). Compact matrix operators between some BK spaces, in: M. Mursaleen (Ed.). {\it{Modern Methods of Analysis and Its Applications, Anamaya Publ., New Delhi}}. 86-120.

\bibitem{article} E. Malkowsky and E. Savas (2004). Matrix transformations between sequence spaces of generalized weighted means. {\it{Applied Mathematics and Computations}}. 147:333-345.

\bibitem{article} F. Basar  and E. Malkowsky (2011). The characterization of compact operators on spaces of strongly summable and bounded sequences. {\it{ Applied Mathematics and Computation}}. 217:5199-5207.

\bibitem{article} G. Darbo (1955). Punti uniti in transformazioni a condominio non compatto. {\it{Rend. Sem. Math. Univ. Padova}}. 24:84-92.

\bibitem{article} H. Kizmaz (1981). On certain sequence spaces. {\it{Canadian Mathematical Bulletin}}. 24(2):169-176.

\bibitem{article} J. Lindenstrauss and L. Tzafriri (1971). On Orlicz sequence spaces. {\it{Israel Journal of Mathematics}}. 10:379-390.

\bibitem{article} K. Kuratowski (1930).  Sur les espaces complets. {\it{ Fund. Math.}}. 15:301-309.

\bibitem{article} M. Mursaleen and A.K. Noman (2010). Compactness by the Hausdorff measure of noncompactness. {\it{ Nonlinear Analysis}}. 73:2541-2557.

\bibitem{article} M. Mursaleen and A.K. Noman (2010). On some new difference sequence spaces of non-absolute type. {\it{Math. Comput. Modelling}}. 52(3-4): 603-617.


\end{thebibliography}
\end{document}